\documentclass[12pt, reqno]{amsart}

\usepackage[utf8]{inputenc}
\usepackage[T1]{fontenc}
\usepackage{tikz,tikz-3dplot}
\usepackage{amsmath,amsfonts,amssymb,amsmath,latexsym,mathrsfs}
\usepackage{shuffle}
\usepackage[all]{xy}
\usepackage{stackengine}
 
\usepackage{enumerate}

\usepackage{hyperref}
 
\usepackage{tabu}
\usepackage{tikz-cd} 

\usepackage{comment}
\usepackage{url}
\usepackage{tikz}

\usepackage{xfrac}

\newtheorem{theorem}{Theorem}[section]

\newtheorem{proposition}[theorem]{Proposition}

\newtheorem{remark}[theorem]{Remark}

\newtheorem{lemma}[theorem]{Lemma}

\newcommand{\Z}{ {\mathbb Z} }
\newcommand{\F}{ {\mathbb F} }
\newcommand{\C}{ {\mathbb C} }

\newcommand{\QZ}{\mathbb{Q}/\mathbb{Z}}

\begin{document}

\title{Second quandle homology of some quandles}

\author{Mohamad \textsc{Maassarani} }
\maketitle
\begin{abstract} We derive formulas from Eisemann's formula for the second quandle homology of quandles. We prove that for a finite injective quandle $Q$, the second quandle homology $H_2(Q)$ can be computed in a representation group of $Q$, which is a finite quotient of the enveloping group of $Q$, by giving a formula for that. We use this to compute $H_2(Q)$ for $Q$ the subquandle of $k$-cycles of $S_n$ and $Q$ the quandle corresponding to the conjugacy class of a diagonalizable element of $SL(n,q)$, with some exceptions. For $G$ a finite group with trivial Bogomolov multiplier, as many classical groups, simple groups and (with exceptions) quasi-simple groups, we derive a formula for $H_2(Conj(G))$ that can be computed in $G$ alone. The formula relies on conjugacy classes of $G$, centralizers in $G$ of representative of such classes and the derived group of $G$. If $G$ is moreover perfect this formula can be related to the first Hochschild homology of the group ring $\Z[G]$. As an application, we compute $H_2(Q)$ for $Q$ the conjugacy quandle of $GL(2,q)$ for $q\neq 2,3$.
\end{abstract}
\section*{Introduction and main results}
 A quandle $Q$ is a set equipped with a binary operation $Q\times Q \to Q$ satisfying some axioms. They originally appeared in knot theory (\cite{Jo}). To a quandle one can associate an enveloping group $G(Q)$ which is infinite for $Q$ non-empty. For quandles one has quandle homology and rack homology. Here, we consider quandle homology, more specifically the second quandle homology group.  Study of second quandle homology is present in the literature. We cite \cite{IV}, \cite{latin},\cite{KI} and \cite{lebed} where the authors consider quandles of transpositions in $S_n$, Alexander quandles, generalized alexander quandles, the conjugacy quandle of the symmetric group and other examples. More generally, Eisermann (\cite{Eis}), gives a formula to compute the second quandle homology from the enveloping group. In this text we derive formulas from this formula.\\\\
A representation group $R(Q)$ for $Q$ is a type of finite quotients of $G(Q)$ (\cite{Ir}). We derive from Eisermann's formula that for $Q$ a connected injective quandle, the second quandle homology is given by : 
$$ H_2(Q)=([R(Q),R(Q)]\cap C_{R(Q)}(q))^{ab},$$  
where $q$ is an element in $R(Q)$ coming from $Q$ and $C_{R(Q)}(q)$ is the centralizer of $q$ in $R(Q)$. This formula is derived in section 1, where we also give some basic definitions (quandles, enveloping groups, representation groups etc.) and derive some facts used in sections $2,3$ and $4$.\\\\
The set of $k$-cycles of $S_n$ forms a quandle with the operation $x\triangleright y=xyx^{-1}$. It is a subquandle of the conjugacy quandle of $S_n$. Such quandles are connected except when $n=k,k+1$ for $k$ odd. We determine representation groups for these injective quandles when they are connected (with the exception of $n=6,7$ for $k$ odd) and compute their second quandle homology using the above formula. The enveloping group for $k$ odd is determined. This is done in section 2, where we treat the cases $k$ even and $k$ odd separately.   \\\\
In section 3, we show that for a finite injective quandle $Q$, the second quandle homology is given by :
$$H_2(Q)=\Z^{n_Q(n_Q-1)}\bigoplus (\bigoplus_i \frac{[{R(Q)},{R(Q)}]\cap C_{R(Q)}({x_i})}{[C_{R(Q)}({x_i}),C_{R(Q)}({x_i})]},$$
where $n_Q$ is the number of orbits of $Q$ and $x_i$ reprensent in $R(Q)$ the orbits of $Q$. We use this formula to compute the second quandle homology for the non connected quandles of $k$-cycles still with the exception of $n=6,7$. We also determine the enveloping groups of these quandles. We note that the last formula for connected differ in form from the one given before. We give  a direct proof that both formulas agree.\\\\
In section 4, we consider $G$ a Schur cover of a centerless group. A conjugacy class of $G$ is a quandle with respect to the binary operation $g\triangleright h=ghg^{-1}$. We show that if the conjugacy class of $x\in G$ ($x$ not in the center) generates $G$, then the conjugacy class form a connected injective quandle, $G$ is a representation group for this quandle and the first formula given here for quandle homology apply. As an application, we determine the second homology of the quandle defined by the conjugacy class of $x\in SL(n,q)$ for $x$ diagonalizable, with some exceptions on $n,q$. \\\\
In section $5$, we show that if a finite group $G$ has trivial Bogomolov multiplier, then the second quandle homology of the conjugacy quandle $Conj(G)$ of $G$ is given by :
$$H_2(Conj(G))\simeq \Z^{c(G)(c(G)-1)} \bigoplus (\bigoplus_{C\in Cl(G)}  \frac{[G,G]\cap C_G(g_C)}{[C_G(g_C),C_G(g_C)]}),$$
where $Cl(G)$ is the set of conjugacy classes of $G$, $g_C$ is a representative of the class $C$ and $C_G(g_C)$ is the centralizer of $g_C$ in $G$. When $G$ is moreover perfect the sum over $Cl(G)$ correspond to the first Hochschild homology group of the group ring $\Z[G]$. We apply the above formula to compute $H_2(Conj(GL(2,q)))$ for $q\neq 2,3$.
\subsection*{AI tools :}  The author used Claude Sonnet 5.5 during the research process and for proof reading. The author checked the information obtained using Claude.

\section{Quandles, second quandle homology and representation groups}
We recall some facts about quandles, second quandle homology and representation groups. We also prove that for an injective connected quandle, the second quandle homology can be computed in a representation group by giving a formula for that. A criterion on wether a Schur cover of the inner automorphism group of a quandle is a representation group for the quandle is given.
  
\subsection{Quandles}
 A quandle is a set $Q$ equipped with a binary operation $\triangleright$  such that :
\begin{itemize}
\item[-] $x\triangleright x=x$ for $x\in Q$.
\item[-] For all $x,y \in Q$ there exist a unique $z\in Q$, such that $x \triangleright z=y$.
\item[-] $x \triangleright( y\triangleright z)= (x \triangleright y) \triangleright (x \triangleright z)$, for $x,y,z \in Q$. 
\end{itemize}
A subquandle of $Q$ is a subset stable under $\triangleright$. One can construct a quandle from a group $G$ using the binary operation defined by $x\triangleright y= xyx^{-1}$  for $x,y\in G$. The defined quandle is called the conjugacy quandle and is usually denoted by $Conj(G)$. A conjugacy class of $G$ is a subquandle of $Conj(G)$. A map $f :Q\to Q'$ is a quandle morphism if $f(x\triangleright y)=f(x)\triangleright f(y)$. It is called isomorphism if $f$ is a bijection.
\subsection{Inner automorphisms}
For $Q$ a quandle and $x\in Q$ denote by $L_x$ the bijection given by $L_x(y)=x\triangleright y$ for $y\in Q$. $L_x$ is called the left translation by $x$. The inner automorphism group $Inn(Q)$ of $Q$ is the subgroup of bijections of $Q$ generated by the bijection $L_x$ for $x\in Q$. The orbits of $Q$ are the orbits of the action of $Inn(Q)$. The morphism $x\mapsto L_x$ is a quandle morphism. A quandle is called faithful if the last map is injective. A quandle is called connected if $Q$ consist of one orbit under the natural action of $Inn(Q)$. For $Q=Conj(G)$ the group $Inn(Q)$ correspond to $Inn(G)$ and the orbits are the conjugacy classes of $G$.
\begin{proposition}\label{In}
Let $G$ be a group and $Q$ be a conjugacy class of $G$. If $Q$ generates $G$ as a group, then $Inn(Q)\simeq Inn(G)$ and $Q$ is connected. The isomorphism is given by $L_x\mapsto c_x$ where $c_x$ is the conjugacy by $x \in G$.  
\end{proposition}
\begin{proof}
$Q$ is stable by $Inn(G)$ and hence we have a morphism $\phi$ from $Inn(G)$ to the bijections of $Q$. Since $Q$ generates $G$, any element of $G$ is a product of elements of $Q$. It follows that any element of $Inn(G)$ restricted to $Q$ is a product of left translations of $Q$. Hence, $\phi$ is in fact a morphism $Inn(G) \to Inn(Q)$. $\phi$ is clearly onto. Now, let $c_g$ be the conjugacy by $g\in G$. Since $Q$ generates $G$, $c_g$ restricts trivially to $Q$ if and only if $g$ is central in $G$. This proves that $\phi$ is injective. We have proved that $\phi: Inn(G)\to Inn(Q)$ is an isomorphism. $Q$ is a single orbit under $Inn(G)$. This proves that $Q$ is connected.  \end{proof}
\subsection{Enveloping group}
The enveloping group $G(Q)$ of a quandle $Q$ is the group :
$$G(Q)=\langle e_x, x \in Q \vert e_xe_ye_x^{-1}=e_{x\triangleright y }, \text{ for } x,y\in Q\rangle.$$
The enveloping group is also known as : the associated group, adjoint group or structure group. The group $G(Q)$ is infinite.
\begin{proposition}
 \item[1)] The map $\varphi_Q : Q \to G(Q)$ mapping $x\in Q$ to the corresponding generator is a quandle morphism for the conjugacy quandle structure on $G(Q)$.
\item[2)]For $f:Q\to Conj(G)$ a quandle morphism there is a unique group morphism $\tilde{f}:G(Q)\to G$ such that $f=\tilde{f}\circ\varphi_Q$. 
\end{proposition}
A quandle $Q$ is called injective if $\varphi_Q : Q \to G(Q)$ is injective.
\begin{proposition}\label{inj}
Let $Q$ be a subquandle of $Conj(G)$ for some group $G$. The quandle $Q$ is injective. In particular, the quandle defined by a conjugacy class of $G$ is injective.
\end{proposition}
\begin{proof}
This follows from the universal property of $\varphi_Q$ applied to the inclusion $Q \to Conj(G)$.
\end{proof}
\begin{proposition}\label{ab}
The abelianisation $G(Q)^{ab}$ of $G(Q)$ is isomorphic to the free abelian group on the set of orbits of $Q$. The abelianisation map can be identified to the morphism mapping $e_x$ to the orbit of $x$.
\end{proposition}
As we have mentioned $Q\to Inn(Q),x\mapsto L_x$ is a quandle morphism. Hence, we have a unique group morphism $\theta_Q :G(Q) \to Inn(Q)$ given by $e_x\mapsto L_x$ for $x\in Q$.
\begin{proposition}[\cite{JD}]
If $Q$ is injective then the kernel of $\theta_Q : G(Q)\to Inn(Q)$ is the center of $G(Q)$.
\end{proposition}
We note that if $Q$ is not injective the kernel of $\theta_Q$ is still central (\cite{Eis}).
\begin{proposition}\label{tor}
Let $Q$ be an injective finite quandle.
\item[1)] The derived group $[G(Q),G(Q)]$ is finite.
\item[2)] Torsion elements of $G(Q)$ form a subgroup, it is $[G(Q),G(Q)]$. 
\end{proposition}
\begin{proof}
By the previous proposition the center of $G(Q)$ has finite index in $G(Q)$. This implies that $[G(Q),G(Q)]$ is finite (Schur theorem for the derived group). We have proved $1)$. The abelianisation of $G(Q)$ is free abelian and hence every torsion element is mapped to zero in the abelianisation. Hence, the torsion elements lie in the derived group and by $1)$ the derived group is finite. This proves $2$.  
\end{proof}
\subsection{Second quandle Homology}
The degree map $d:G(Q) \to Z$ is the group morphism mapping $e_x \to 1$ for $x\in Q$.
Eisermann gives a formula for computing the second homology of any quandle from the enveloping groups (\cite{Eis}, theorem 9.9). For injective quandle we get : 
\begin{theorem}[\cite{Eis}]\label{eis0}
Let $Q$ be an injective quandle. For each orbit of $Q$ choose an element $x_i $ of the orbit. The second quandle homology of $Q$ is given by : 
\[H_2(Q)\simeq \oplus_i (C_{G(Q)}(e_{x_i})\cap \mathrm{ker(d)})^{ab},\]
where $C_{G(Q)}(e_{x_i})$ is the centralizer of $e_{x_i}$ in $G(Q)$.
\end{theorem}
\begin{proposition}
For $Q$ connected the kernel of the degree map is $[G(Q),G(Q)]$.
\end{proposition}
\begin{proof}
The degree map can be identified to the abelianisation map for $Q$ connected.
\end{proof}
Combining the last theorem and the last proposition, we get.
\begin{proposition}\label{eis}
Let $Q$ be a connected injective quandle. For any $x\in Q$, the second quandle homology of $Q$ is given by :
$$H_2(Q)\simeq (C_{G(Q)}(e_x)\cap [G(Q),G(Q)])^{ab}.$$
where $C_{G(Q)}(e_x)$ denotes the centralizer of $e_x$ in $G(Q)$ and $x$ is any element of $Q$.
\end{proposition}

 \subsection{Stem extensions}
The Schur multiplier $M(G)$ of a finite group $G$ is the cohomology group $H^2(G,\C^\times)$. The Schur multiplier is finite. A stem extension is a central extension $1 \to Z \to E \to G \to 1$ where the nucleus $Z$ is a subgroup of $ [E,E]$. The inflation-restriction five-term exact sequence for a stem extension $1\to Z \to E \to G\to 1$ gives an exact sequence :
$$0 \to Hom(Z,\mathbb{C}^\times)\to H^2(G,\mathbb{C}^\times)\to H^2(E,\mathbb{C}^\times).$$
In particular $Z$ is finite and $Z$ correspond to a subgroup of $M(G)$. Finally, a Schur cover is a stem extension where $Z\simeq M(G)$.

\subsection{Representation groups}\label{repgr}
In this subsection, we consider a finite injective quandle $Q$. We will recall the definition of a representation group from \cite{Ir}. The notion is defined for non injective finite quandles also. We have seen that the center $Z_Q$ of $G(Q)$ is the kernel of $\theta_Q : G(Q) \to Inn(Q)$. Hence, $Z_Q$ is of finite index in $G(Q)$ and hence $Z_Q$ is finitely presented since $G(Q)$ is finitely presented. This proves that $Z_Q$ is a finitely generated abelian group. We have seen that the derived group of $G(Q)$ is finite and consist of the torsion elements of $G(Q)$ (proposition \ref{tor}) implies that $[G(Q),G(Q)]\cap Z_Q=Tor(Z_Q)$ where $Tor(Z_Q)$ is the torsion subgroup of $Z_Q$. Knowing this and that $Z_Q$ has finite index in $G(Q)$ (hence the image of $Z_Q$ has finite index in $G(Q)^{ab}$), we get that $Z_Q \simeq \Z^{n_Q} \bigoplus Tor(Z_Q)$, where $n_Q$ is the number of orbits of $Q$ under $Inn(Q)$.\\\\
Let $\alpha : \Z^{n_Q} \bigoplus Tor(Z_Q) \to Z_Q$ be an isomorphism. The representation group of $Q$ associated $\alpha$ is the group $G(Q)_\alpha=G(Q)/\alpha(\Z^{n_Q})$. We have a unique morphism $\theta_{Q,\alpha} : G(Q)_\alpha \to Inn(Q)$ such that $\theta_{Q,\alpha}\circ \pi_\alpha=\theta_Q$, where $\pi_\alpha$ is the quotient map $G(Q)\to G(Q)_\alpha$. We note that $\theta_{Q,\alpha}$ defines a stem extension of $Inn(Q)$ with nucleus isomorphic to $Tor(Z_Q)$, that $G(Q)_\alpha$ is finite and that $\pi_\alpha$ restricts to an isomorphism $[G(Q),G(Q)]\to [G(Q)_\alpha,G(Q)_\alpha]$ since $\mathrm{ker}(\pi_\alpha) \cap [G(Q),G(Q)]=1$. In our work, we will say that $(G,f)$, where $f:G\to Inn(Q)$ is a stem extension, is  a representation group if we have an isomorphism $g :G \to  G(Q)_\alpha$ for some $\alpha$ such that $ \theta_{Q,\alpha}\circ g=f$. 
\begin{lemma}\label{tor1}
For $Q$ and injective finite quandle, $Tor(Z_Q)$ can be identified to a subgroup of $M(Inn(Q))$.
\end{lemma}
\begin{proof}
$G(Q)_\alpha$ is a stem extension of $Inn(Q)$. Hence, the nucleus of the extension $G(Q)_\alpha \to Inn(Q)$ correspond to a subgroup of $M(Inn(Q))$. The lemma follows from the fact that the nucleus of $G(Q)_\alpha\to Inn(Q)$ is isomorphic to $Tor(Z_Q)$. 
\end{proof}
\begin{lemma}\label{Sub}
Let $G\to H$ be a stem extension and let $K$ be a subgroup of $G$ mapping onto $H$. We have $K=G$.
\end{lemma}
\begin{proof}
Denote by $Z$ the nucleus of the stem extension. We have $G=KZ$. Since $Z$ is central $[G,G]=[K,K]$. Since $Z\subset [G,G]$, $Z\subset K$. This proves that $K=G$, since $G=KZ$.
\end{proof}
\begin{proposition}\label{RS}
Let $Q$ be a finite injective quandle and $f: G\to Inn(Q)$ be a Schur cover. $(G,f)$ is a representation group for $Q$ if and only if there exist a morphism $g : G(Q) \to G$ such that $f\circ g = \theta_Q$. 
\end{proposition}
\begin{proof}
The implication "$G$ is a representation group implies the existence of $g$" is clear. We prove the converse. Let $g: G(Q)\to G$ be as in the proposition. The image of $g$ maps under $f$ onto $Inn(Q)$ hence by the previous lemma $g$ is surjective. $f\circ g(Tor(Z_Q))=\theta_Q(Tor(Z_Q))=1$ hence $g$ maps $Tor(Z_Q)$ to $Z$ the nucleus of $f$. It follows from the surjectivity of $g$ and the fact that $f$ is a stem extension that we have the following sequence of morphisms induced by $f$ and $g$ :
$$ G(Q)'/Tor(Z_Q) \to G'/Z \to Inn(Q)',$$
where $H'$ denotes $[H,H]$. The last map is an isomorphism and the composite is an isomorphism since $Tor(Z_Q)=Z_Q\cap G(Q)'$. Hence, $G(Q)'/Tor(Z_Q) \to G'/Z$ is an isomorphism. This proves that the restriction of $g$, $Tor(Z_Q)\to Z$ is surjective. Now since $f$ is a Schur cover $Z\simeq M(Inn(Q))$. Since $Tor(Z_Q)$ identifies to a subgroup of $M(Inn(Q))$ (lemma \ref{tor1}), the last two remarks imply that $g$ maps $Tor(Z_Q)$ onto $Z$ isomorphically. Let $K$ be the kernel of $g$. Since $f\circ g =\theta_Q$, $K\subset Z_Q$. Since $g$ maps $Tor(Z_Q)$ isomorphically to $Z$. $K$ contains no torsion and is hence free abelian. The group $L=\langle K,Tor(Z_Q) \rangle\subset Z_Q$ is isomorphic to $K\oplus Tor(Z_Q)$. Now we have a sequence of morphisms induced by $g$ and $f$ :
$$G(Q)/L\to G/Z \to Inn(Q).$$
$G(Q)/L$ to $G/Z$ is an isomorphism and so is $G/Z \to Inn(Q)$. Hence, $G(Q)/L\to Inn(Q)$ is an isomorphism proving that $L\simeq Z_Q$. Hence, we have an isomorphism $\beta :\Z^{n_Q}\to K$ let $\alpha : \Z^{n_Q}\oplus Tor(Z_Q) \to L=Z_Q$ be the isomorphism extending $\beta$ by the identity isomorphism on $Tor(Z_Q)$. We have naturally an isomorphism $h$ between $G(Q)_\alpha=G(Q)/K$ and $G$, and $\theta_{Q,\alpha}=f\circ h$. This proves that $(G,f)$ is a representation group. 
\end{proof}
\begin{proposition}\label{RH}
Let $(G,f)$ be a representation group for $Q$ a finite injective and connected quandle. For any $x \in Q$, the second quandle homology of $Q$ is given by :
$$H_2(Q)\simeq (C_{G}(y)\cap [G,G])^{ab}.$$
for any $y\in f^{-1}(L_x)$ with $C_{G}(y)$ the centralizer of $y$ in $G$.
\end{proposition}
\begin{proof}
We prove the proposition for the representation groups of the form $(G(Q)_\alpha,\theta_{Q,\alpha})$, since the proposition is consequence of this. We fix $\alpha$. Take $x \in Q$. We have by proposition \ref{eis} :
$$H_2(Q)\simeq (C_{G(Q)}(e_x)\cap [G(Q),G(Q)])^{ab}.$$
for any $x\in Q$. We have denoted by $\pi_\alpha $ the projection $G(Q)\to G(Q)_\alpha$. Denote by $K$ the kernel of $\pi_\alpha$. We have see in subsection \ref{repgr}, that  $\pi_\alpha $ resticts to an isomorphism $G(Q)'\to G(Q)_\alpha'$, where $H'=[H,H]$. Hence we have an injective homomorphism $$g: C_{G(Q)}(e_x)\cap G(Q)'\to C_{G(Q)_\alpha}(\pi_\alpha(e_x))\cap G(Q)_\alpha',$$ induced by $\pi_\alpha$. Suppose that $g$ is not surjective. This means there is $z\in C_{G(Q)_\alpha}(\pi_\alpha(e_x))\cap G(Q)_\alpha'$. Let $\tilde{z}$ be a lift of $z$ to $G(Q)'$. We have $\pi([\tilde z, e_x])=[z,\pi_\alpha(e_x)]=1$. This implies that $[\tilde{z},e_x]\in K\cap G(Q)'=1$. This proves that $\tilde z$ also lies in $C_{G(Q)}(e_x)$ and hence $z$ is in the image of $g$. $g$ is therefore surjective. We proved already that $g$ is injective. We have that $g$ is an isomorphism. It follows that for any $x\in Q$ :
$$H_2(Q)\simeq (C_{G(Q)_\alpha}(\pi_\alpha(e_x))\cap [G(Q)_\alpha,G(Q)_\alpha])^{ab}.$$
This proves the proposition for $G(Q)_\alpha$ for a given lift of $L_x$ for any $x$. But lifts of $L_x$ differ by a central elements and hence define the same group.
\end{proof}
\section{connected $k$-cycle quandles}
We consider within the subquandles of $Conj(S_n)$ consisting of $k$-cycles the connected ones. We deterrmine representation groups and compute the second quandle homology for these quandles with the exceptions of $n=6,7$ when $k$ is odd. We also determine the enveloping group for $k$ odd. The first subsection is dedicated for Schur covers of $S_n$ and $A_n$. We treat the case $k$ even in the second subsection then the case $k$ odd in the third subsection. 

\subsection{Schur covers of symmetric and alternating groups}
We will recall some elements related to Schur covers of $S_n$ and $A_n$. They are due essentially to Schur, and they can be found in \cite{Sc}. We will also prove some facts in the covers.

\begin{theorem}[Schur]
\item[1)] For $n\geq 3$, the Schur multiplier of $S_n$ is given by :
\[M(S_n) =
\begin{cases}
  1 & \text{if } n = 3,\\
    
  \Z/2\Z & \text{otherwise.}
\end{cases}\]
\item[2)] For $n\geq 3$, the Schur multiplier of $A_n$ is given by :
\[M(A_n) =
\begin{cases}
  1 & \text{if } n = 3,\\
    \Z/6\Z & \text{if } n=6,7, \\
  \Z/2\Z & \text{otherwise.}
\end{cases}\]
\end{theorem}
For $n\geq 4$. Let $2\cdot S_n^+$ and $2\cdot S_n^-$ be the groups definied by :

$$ \langle z,t_1,\dots,t_{n-1} \vert z^2=1,\: zt_i=t_iz,\: t_i^2=z^\epsilon,\: (t_it_{i+1})^3=z^\epsilon, \:(t_it_j)^2=z \: (\vert i-j \vert \geq 2) \rangle,$$
where $\epsilon=1$ for $2\cdot S_n^{+}$ and $\epsilon=0$ for $2\cdot S_n^-$. We have morphisms $p_+: 2 \cdot S_n^+\to S_n$ and $p_-:2\cdot S_n^- \to S_n$ given by :
$$p_+(z)=p_-(z)=1,\: p_+(t_i)=p_-(t_i)=(i,i+1).$$
The kernel of both maps is the order $2$ subgroup $\langle z \rangle$, it correspond to the center of these groups. 
\begin{theorem}[Schur]
For $n\geq 4$, $S_n$ has exactly two Schur covers up to isomorphism of extensions of $S_n$, they are $(2\cdot S_n^+,p_+)$ and $(2\cdot S_n^-,p_-)$.
\end{theorem}
For $n=6$, the groups in both covers of $S_n$ are isomorphic (Schur) but the covering maps are not. For any $n$, lifts of $(i,i+1)$ have order $4$ in $2\cdot S_n^{+}$, while they have order $2$ in $2\cdot S_n^-$.\\\\
Denote by $2\cdot A_n$ the preimage of $A_n$ in $2\cdot S_n^+$ ($n\geq 4$). $p_+$ resticts to a morphism $2\cdot A_n \to A_n$. We will also denote it $p_+$.
\begin{proposition}[Schur]
For $n\geq 4$ and $n\neq 6,7$, up to isomorphism $A_n$ has a unique Schur cover it is given by $(2\cdot A_n,p_+)$.
\end{proposition}

Recall that the parity $p(\sigma)$ of $\sigma \in S_n$ is $0$ if the sign of $\sigma$ is $1$ and it is equal to $1$ if the sign of $\sigma$ is $-1$. A permutation is odd if its parity is one and it is even if the parity is $0$. We will use the commutator notation $[a,b]=aba^{-1}b^{-1}$.

\begin{proposition}\label{com}
Let $n\ge4$ and $S$ be $2\cdot S_n^{+}$ or $2\cdot S_n^-$. Let $\sigma,\tau\in S_n$ be permutations with disjoint supports and parities $p(\sigma),p(\tau)\in\{0,1\}$. For any lifts $\tilde\sigma,\tilde\tau$ to $S$,
\[
[\tilde\sigma,\tilde\tau]=z^{\,p(\sigma)p(\tau)}.
\]
 
\end{proposition}
\begin{proof}
For $j,k$ with $|j-k|\ge2$, we have $[t_j,t_k]=t_jt_kt_j^{-1}t_k^{-1}=(t_jt_k)^2=z$, using $t_j^{-1}=z^{\epsilon}t_j$ and $z^2=1$. Any two transpositions with disjoint supportsf are conjugate in $S_n$ to $((12),(34))$; conjugating by a lift of the conjugating permutation fixes $z$ and carries commutators of lifts to commutators of lifts. Changing a lift multiplies by a central element and does not change the commutator. So $[\tilde u,\tilde v]=z$ for disjoint transpositions $u,v$. Now let $u$ and $v$ be as in the proposition. Write $\sigma=u_1\cdots u_r$ and $\tau=v_1\cdots v_s$ as products of transpositions with $r\equiv p(\sigma)$ mod $2$, $s\equiv p(\tau)$ mod $2$, all $u_i$ supported in $\operatorname{supp}\sigma$ and all $v_j$ in $\operatorname{supp}\tau$. Each $[\tilde u_i,\tilde v_j]=z$ is central, so $$[ab,c]=a[b,c]a^{-1}[a,c]\quad \text{and} \quad [a,bc]=[a,b]b[a,c]b^{-1},$$ give that $[\tilde\sigma,\tilde\tau]=z^{rs}$. 
\end{proof}
\begin{proposition}
For $n=6,7$, $p_+: 2\cdot A_n \to A_n$ is a stem extension.
\end{proposition}
\begin{proof}
Let $u_1,u_2,u_3$ be lift of the transposition $(12),(34),(56)$ respectively. $u_1u_2$ and $u_2u_3$ lie in $2\cdot A_n$.
Using the identities : $$[ab,c]=a[b,c]a^{-1}[a,c]\quad \text{and} \quad [a,bc]=[a,b]b[a,c]b^{-1},$$
we get that $$[u_1u_2,u_2u_3]=u_1[u_2,u_2u_3]u_1^{-1}[u_1,u_2u_3]=u_1[u_2,u_2]u_2[u_2,u_3]u_2^{-1}u_1^{-1}[u_1,u_2]u_2[u_1,u_3]u_2^{-1}.$$
By the previous proposition $[u_i,u_j]=z$ for $i\neq j$. Hence, 
$$[u_1u_2,u_2u_3]=z.$$
This proves that $z\in [2\cdot A_n,2\cdot A_n]$. We have proved the proposition.
\end{proof}
\begin{lemma}
Let $x$ be a $k$-cycle with $k$ even and let $\tilde x$ be a lift of $x$ to $2\cdot S_n^-$. We have :$$\tilde x^k=z^{\frac{k(k-2)}{8}}.$$
\end{lemma}
\begin{proof}
$x^{k/2}$ is a product of $l=k/2$ disjoint transposition. Let $u_1,\dots,u_l$ be lifts of these transposition. We have that either $\tilde x^{l}=u_1\cdots u_l$ either $\tilde x^l=zu_1\cdots u_l$. In both cases we have :
$$(\tilde x^l)^2=\tilde x^k=(u_1\cdots u_l)^2=u_1\cdots u_l u_1\cdots u_l.$$
By the previous proposition $[u_i,u_j]=z$ for $i\neq j$ and hence $u_iu_j=zu_ju_i$. The relation $t_i^2=1$ implies that $u_j^2=1$, since up to multiplication by $z$, $t_i$ and $u_j$ are conjugate. Combining the relations $u_iu_j=zu_ju_i$ and $u_j^2$, we can reduce the last equation to :
$$\tilde x^k=z^{\frac{l(l-1)}{2}}=z^{\frac{k(k-2)}{8}}.$$
\end{proof}
\begin{proposition}\label{x2}
Let $x$ be a $k$-cycle with $k$ even and $\tilde x$ be a lift of $x$ to $2\cdot S_n^-$. The subgroup generated by $\tilde x^2$ and $z$ is given by :
$$\langle \tilde x^2,z \rangle \simeq \begin{cases}
  \Z/(k/2)\Z \times \Z/2\Z & \text{if }k\equiv0\pmod 8,\\
    
  \Z/k\Z & \text{otherwise.}
\end{cases}$$
\end{proposition}
\begin{proof}
The group $\langle \tilde x^2,z \rangle$ is abelian. If $\frac{k(k-2)}{8}$ is odd, by the previous lemma. $(\tilde x^2)^{k/2}=z$. Hence, $\tilde x^2$ has order $k$ and $\langle \tilde x^2,z \rangle=\langle \tilde x^2 \rangle\simeq \Z/k\Z$. If $\frac{k(k-2)}{8}$ is even, we have, by using the previous lemma, that $\tilde x^2$ has order $k/2$ and $\langle \tilde x^2 \rangle$ maps isomorphically onto $\langle x^2 \rangle$ and hence $z\notin \langle \tilde x^2 \rangle$. This proves that, for $\frac{k(k-2)}{8}$ even, $\langle \tilde x^2,z \rangle \simeq \Z/(k/2)\Z \times \Z/2\Z $. For $k$ congruent to $4$ or $6$ mod $8$, $\frac{k(k-2)}{8}$ is odd an hence $\langle \tilde x^2,z \rangle\simeq \Z/k\Z$. If $k$ is congruent to $0$ or $2$ mod $8$, then $\frac{k(k-2)}{8}$ is even and hence $\langle \tilde x^2,z \rangle \simeq \Z/(k/2)\Z \times \Z/2\Z $. Finally, for $k$ congruent to $2$ mod $8$, $k/2$ is odd and hence  $\Z/(k/2)\Z \times \Z/2\Z \simeq \Z/k\Z$.
\end{proof}

\subsection{$k$-cycles for k even}
In this subsection, $k$ and $n$ are integers with $2\leq k\leq n$, $n\geq 3$ and $k$ even. $Q$ will denote the subquandle of the conjugacy quandle of $S_n$ consisting of $k$-cycles. $k$-cycles for $k$ even generate $S_n$ and form a conjugacy class. Hence, by proposition \ref{In}, $Q$ is connected and $Inn(Q)\simeq Inn(S_n) \simeq S_n$. We will identify $Q\to Inn(Q),x\mapsto L_x$ to the inclusion $Q\to S_n$. The quandle $Q$ is also injective by proposition \ref{inj}.

\begin{theorem}
\item[1)] If $n-k\geq 2$ or $(n,k)=(3,2)$, then $S_n$ is the representation group of $Q$ the quandle of $k$-cycles in $S_n$ with $k$ even.
\item[2)] If $n-k=0,1$ with $(n,k)\neq (3,2)$, then $2\cdot S_n^+$ and $2 \cdot S_n^-$ are representation groups of $Q$ the quandle of $k$-cycles in $S_n$ with $k$ even.
\end{theorem}
\begin{proof}
We prove $1)$. For $k=2$ and $n=3$, $M(S_n)=0$. Since a representation group is a stem extension of $Inn(Q)\simeq S_n$, the only possible representation group is $S_n$. Now take, $n\geq 4$. $M(S_n)=\Z/2\Z$. The possible representation groups for $Q$ are hence up to isomorphism $2\cdot S_n^+$, $2\cdot S_n^-$ or $S_n$. We will prove that there is no morphism from $G(Q)$ to any double cover of $S_n$ mapping $e_x$ to a lift of $x$. This proves, taking into consideration proposition \ref{RS}, that $S_n$ is the representation group of $Q$. Now suppose that we have a morphism from $G(Q)$ to one of the double covers mapping $e_x$ to a lift $\tilde{x}$ of $x$. Put $\tilde{C}=\{\tilde{x}\vert x\in Q\}$. By lemma \ref{Sub} the morphism $G(Q)$ to the double cover is surjective. Since conjugacy in $G(Q)$ permutes the generators. Conjugacy in the double cover permutes $\tilde{C}$. The projection from the double cover to $S_n$ restricts to a bijection from $\tilde{C}$ to $Q$. Now fix $x\in Q$ and let $g$ be a transposition disjoint for $x$ ($n-k \geq 2$). $[g,x]=1$. Let $\tilde{g}$ be a lift of $g$. $\tilde{g}\tilde{x}\tilde{g}^{-1}$ lies in $\tilde{C}$ and lies over $gxg^{-1}=x$. Hence, $\tilde{g}\tilde{x}\tilde{g}^{-1}=\tilde{x}$ and $[\tilde{g},\tilde{x}]=1$. This contradicts proposition \ref{com}, since $p(x)=p(g)=1$. We have proved $1)$. \\
We prove $2)$. Take $n-k=0,1$, fix $x\in Q$ and let $S$ be any of the double covers of $S_n$ ($n$ is greater then $4$ since $(n,k)\neq (3,2)$). The centralizer $C_{S_n}(x)$ is equal to $\langle x\rangle$. Let $\tilde x$ be a lift of $x$ to $S$. The preimage of $\langle x\rangle$ is generated by $\tilde x$ and the central $z$, hence abelian, so $C_{S}(\tilde x)$ is the whole preimage of $C_{S_n}(x)$, of order $2k$. The class of $\tilde x$ has size $2\cdot (n!)/(2k)=n!/k=\vert Q\vert$, so it maps bijectively onto $Q$ and is a quandle isomorphic to $Q$. Hence, we have a morphism from $G(Q) \to S$ which composition with the projection to $S_n$ is given by $e_y\mapsto y$. It follows from proposition \ref{RS} that $S$ is a representation group of $Q$.
\end{proof}

\begin{theorem}
The second quandle homology of the quandle $Q$ of $k$-cycles in $S_n$ with $k$ even, is given by:
\item[1)] For $(n,k)=(3,2)$ : $$H_2(Q)=0.$$
\item[2)] For $n-k\geq 2$:  $$H_2(Q)=\Z/k\Z.$$
\item[3)]For $n-k=0,1$, $(n,k)\neq (3,2)$ :
$$H_2(Q)=\begin{cases}
  \Z/(k/2)\Z \times \Z/2\Z & \text{if }k\equiv0\pmod 8,\\
    
  \Z/k\Z & \text{otherwise.}
\end{cases}.$$
\end{theorem}
\begin{proof}
We apply proposition \ref{RH} to compute the second quandle homology from the representation group. We prove $1)$. By the previous theorem the representation group for $Q$ is $S_3$ for $(n,k)=(3,2)$. We should consider the intersection $(C_{S_3}(x)\cap [S_3,S_3])^{ab}$ for $x$ a transposition in $S_3$. The centralizer consists of the subgroup generated by the transposition and hence intersect $A_3$ trivially. This proves that $H_2(Q)=0$, for $(n,k)=(3,2)$.\\
We now prove $2)$. For $n-k\geq 2$, the representation group of $Q$ is $S_n$ and the second quandle homology is $K^{ab}$ for $K=C_{S_n}(x)\cap A_n$ and $x$ any $k$-cycle. $C_{S_n}(x)\simeq \langle x\rangle \times S_{n-k}$ and its intersection with $A_n$ correspond to the elements $x^iy$ where $y$ are elements in the $S_{n-k}$ factor, with $i \equiv p(y) \pmod 2$ where $p(y)$ is the parity of $y$. Now the commutator of $x^iy$ and $x^jy'$ is given by $[y,y']$ this proves that $[K,K]\simeq 0 \times A_{n-k}$. Hence, $K^{ab}$ is the subgroup of $\langle x \rangle \times S_{n-k}/\{1\}\times A_{n-k} \simeq \Z/k\Z\times \Z/2\Z$ consisting of pairs $(l,r)$ with $l \equiv r \pmod 2$. The projection morphism $\Z/k\Z\times Z/2\Z \to \Z/k\Z$ restricts to an isomorphism between $K^{ab}$ and $\Z/k\Z$.\\
We will prove $3)$. For $n-k=0,1$, $(n,k)\neq (3,2)$, $2 \cdot S_n^-$ is a representation group for $Q$. We can consider the abelianisation of $K=C_{2\cdot S_n^-}(\tilde x)\cap [ 2\cdot S_n^-,2\cdot S_n^-]$ for $\tilde x $ a lift of a $k$-cycle $x$. We have seen in the proof of $2)$ of the previous theorem that $C_{2\cdot S_n^-}(\tilde x)=\langle \tilde x,z\rangle$. Now since $  2\cdot S_n^-$ is  a stem extension $[ 2\cdot S_n^- ,2\cdot S_n^-]=p_-^{-1}(A_n)$. Hence, $K=C_{2\cdot S_n^-}(\tilde x)\cap [ 2\cdot S_n^-,2\cdot S_n^-]=\langle \tilde x^2,z\rangle$. $K$ is abelian and hence $H_2(Q)\simeq K^{ab}=\langle \tilde x^2,z\rangle$. The groups in $3)$ correspond to the computations of $\langle \tilde x^2,z\rangle$ given by proposition \ref{x2}.

\end{proof}

\subsection{$k$-cycles for k odd}
In this subsection, $k$ and $n$ are integers with $0<k\leq n$, $n\geq 3$ and $k\geq 2$ is odd. We also assume that $n\neq 6,7,k, k+1$. This implies that $n\geq 5$. We will denote by $Q$ the subquandle of $S_n$ consisting of $k$-cycles. It is also a subquandle of $A_n$. The set $Q$ is a conjugacy class of $A_n$ generating $A_n$ and hence, by proposition \ref{In}, $Q$ is connected and $Inn(Q)\simeq Inn(A_n)\simeq A_n$. We will identify the map $Q\to Inn(Q),x\mapsto L_x$ to the inclusion of $Q$ in $A_n$. The quandle $Q$ is also injective by proposition \ref{inj}

\begin{theorem}\label{x*}
For $k$ odd and $n\neq 6,7,k,k+1$, the group $2\cdot A_n$ is a representation group of $Q$ the subquandle of $k$-cycles of the conjugacy quandle of $S_n$.
\end{theorem}
\begin{proof}
Take $x\in Q$ and let $\tilde x$ be a lift of $x$ to $2\cdot A_n$. The two lifts of $x$ are $\tilde  x $ and $z\tilde x$. $\tilde x^k \in \{1,z\}$ and $(z\tilde x)^k=z\tilde x^k$. Hence, only one of the lifts of $x$ has order $k$ denote by $x^*$ this lift. Take $y^*$ the unique lift of order $k$ of $y \in Q$. $y^*x^*(y^*)^{-1}$ is a lift over $yxy^{-1}$ and is of order $k$. Hence, $(yxy^{-1})^*=y^*x^*(y^*)^{-1}$. This proves that we have a quandle isomorphism $Q \to C$ where $C=\{x^* \vert x\in Q\}$ given by $x\mapsto x^*$. It follow that we have a morphism $G(Q) \to 2\cdot A_n$ given by $e_x\mapsto x^*$ for $x\in Q$. The composition of this morphism by $p_+$ is given by $e_x\mapsto x$ for $x\in Q$. It follows from proposition \ref{RS} that $2\cdot A_n$ is a representation group for $Q$.
\end{proof}
\begin{remark}
Let $x\neq 1$ be an element of order $k$ of a simple group $G$. The conjugacy class $Cl_G(x)$ of $x$ generates a normal subgroup of $G$ and hence generates the group $G$. $Cl_G(x)$ seen as a subquandle of $Conj(G)$ is hence injective connected and has inner automorphism group isomorphic to $G$. If $k$ and $\vert M(G) \vert$ are comprime then one can prove that $x$ admit a unique lift of order $k$ in the Schur cover of $G$. This can be checked by considering $k$-th powers of the possible lifts as we have done in the previous proof. In particular, by adapting the last proof we get that the Schur cover of $G$ is a representation group of the quandle $Cl_G(x)$. This for instance apply to the quandle $5$-cycles of $A_7$ and for the conjugacy class of any element of odd order in $A_n$ with $n\neq 6,7$.    
\end{remark}
\begin{lemma}
Let $\tilde{C}$ be the lift of the centraliser of a $k$-cycle in $A_n$ to $2\cdot A_n$.

$$ \tilde C^{ab}=\begin{cases}
\Z/k\Z\times \Z/2\Z & \text{if } n-k=2,\\
\Z/k\Z\times \Z/6\Z & \text{if } n-k=3,\\
\Z/k\Z\times \Z/3\Z & \text{if } n-k=4,\\
\Z/k\Z & \text{if } n-k\geq 5
\end{cases}$$
\end{lemma}
\begin{proof}
Let $x$ be a $k$-cycle in $A_n$. The centralizer of $x$ in $S_n$ is isomorphic to $\langle x \rangle \times S_{n-k}$. Since $x$ is even, the centralizer in $A_n$ is isomorphic to $\langle x \rangle \times A_{n-k}$. Denote by $A$ the group corresponding to the factor $A_{n-k}$. Note that we have seen in the last proof that $x$ admit a lift $x^*$ to $2\cdot A_n$ of order $k$. We have $\tilde C=\langle x^* \rangle p_+^{-1}(A)$ (it follows from proposition \ref{com}, that $\langle x^* \rangle$ is central in $\tilde C$). The image of $\langle x^* \rangle \cap p_+^{-1}(A)$ under $p_+$ lies in $\langle x \rangle \cap A=1$. Hence, $\langle x^* \rangle \cap p_+^{-1}(A)\subset \langle z\rangle$. But $\langle x^* \rangle \cap \langle z \rangle =1$ ($x^*$ has order $k$) and therefore  $\langle x^* \rangle \cap p_+^{-1}(A)=1$. This proves that $\tilde C\simeq \langle x^* \rangle \times p_+^{-1}(A)$. For $n-k=2$, $A=1$, $p_+^{-1}(A)\simeq \Z/2\Z$. The case $n-k=2$ in the lemma follows. For $n-k=3$, $A$ is cyclic of order $3$ its lift is of order $6$. Take a $3$-cycle $y$ in $A$ and consider its lift $zy^*$. $zy^*$ is of order $6$ and hence $p_+^{-1}(A) \simeq \Z/6\Z$. It follows that $\tilde C^{ab} \simeq \Z/k\Z\times \Z/6\Z$ for $n-k=3$. Now take $n-k \geq 4$. With no loss of generality, we suppose that $x=(n-k+1, n-k+2,\dots,n)$. Denote by $f: 2\cdot S_{n-k}^+\to 2\cdot S_n^+$ the morphism mapping $t_i$ to $t_i$ for $i\leq n-k-1$ and $z\to z$. Let $S_{n-k}$ be the symmetric group on $\{1,..,n-k\}$ seen as a subgroup of $S_n$. The image of $f$ maps onto $S_{n-k}$ surjectively and the image contains $z$ of order $2$ that is mapped to $0$ in $S_{n-k}$. Hence, $f$ is injective. As one checks, $p_+^{-1}(A)$ (for $p_+:2\cdot S_n^+ \to S_n$, here $A$ is the alternating group on $\{1,..,n-k\}$) is equal to $f(2\cdot A_{n-k})$, since $p_+^{n-k}= p_+^n\circ f$ where $p_+^l$ is the cover map $p_+: 2\cdot S_l^+ \to S_l$. This proves that $p_+^{-1}(A)\simeq 2\cdot A_{n-k}$. Hence, we have $\tilde C  \simeq \Z/ k\Z \times 2\cdot A_{n-k}$ for $n-k\geq 4$. Since $2\cdot A_{n-k}$ is a stem extension of $A_{n-k}$, $(2\cdot A_{n-k})^{ab} \simeq A_{n-k}^{ab}$. For $n-k\geq 5$, $A_{n-k}^{ab}=1$ and for $n-k=4$, $A_{n-k}^{ab}\simeq \Z/3\Z$. This completes the proof of the lemma.
\end{proof}
\begin{theorem}
The second quandle homology of the quandle $Q$ of $k$-cycles in $S_n$ with $k\geq 2$ odd $($$n\neq 6,7,k,k+1$$)$, is given by : 
$$ H_2(Q)=\begin{cases}
\Z/k\Z\times \Z/2\Z & \text{if } n-k=2,\\
\Z/k\Z\times \Z/6\Z & \text{if } n-k=3,\\
\Z/k\Z\times \Z/3\Z & \text{if } n-k=4,\\
\Z/k\Z & \text{if } n-k\geq 5
\end{cases}$$
\end{theorem}
\begin{proof}
$2\cdot A_n$ is a representation group for $Q$. With the conditions on $(n,k)$, we have $n\geq 5$. Hence,  $2\cdot A_n$ is perfect in this range since it is a stem extension of the perfect group $A_n$. Using the last two remarks and applying proposition \ref{RH}, we get :
$$ H_2(Q)=C_{2\cdot A_n}(\tilde x)^{ab},$$
for $\tilde x$ a lift of any element $x\in Q\subset A_n$. Since $2\cdot A_n\to A_n$ is surjective with central kernel $\langle z \rangle$, we get an exact sequence :
$$0\to \langle z \rangle \to C_{2\cdot A_n}(\tilde x) \to C_{A_n}(x).$$
Now for $y \in Q$ denote by $y^*$ the unique lift of $y$ to $2\cdot A_n$ of order $k$ (it exists for $k$ odd as we have seen in previous proofs). Set $A=\{y^* \vert y\in Q \}$ and take $g\in 2\cdot A_n$, $gy^*g^{-1}$ is a lift of order $k$ of $p_+(g)yp_+(g)^{-1}$. Hence, $gy^*g^{-1}$ lies in $A$. This proves that the $A$ is stable by conjugacy in $2\cdot A_n$. Moreover, $Q\to A,y\mapsto y^*$ is a quandle isomorphism. Since $Q$ is connected $A$ is a connected quandle and hence $A$ is a conjugacy class with cardinal $\vert Q\vert$. It follows that :
$$\vert A \vert =\vert 2\cdot A_n \vert / \vert C_{2 \cdot A_n}(x^*)\vert = \vert Q\vert =\vert A_n \vert / \vert C_{A_n}(x)\vert, $$
and hence  $\vert C_{2 \cdot An}(x^*)\vert= 2 \vert C_{A_n}(x) \vert $. Using the last equality and the exact sequence, we see that $C_{2\cdot A_n}(x^*)$ is the preimage of $C_{A_n}(x)$ under $p_+$. Since $H_2(Q)\simeq C_{2\cdot A_n}(x^*)^{ab}$, the theorem follows from the previous lemma.
\end{proof}
\begin{theorem}\label{envCycle}
For $n\geq 5$, the enveloping group of the $k$-cycle quandle $Q$ in $S_n$ is isomorphic to $2 \cdot A_n \times \Z$.
\end{theorem}
\begin{proof}
We have seen $2\cdot A_n$ is a representation group for $Q$. Hence we have a surjective morphism $f :G(Q)\to 2\cdot A_n$ with kernel $K$ central and isomorphic to $\Z$ (the quandle is connected and the number of orbits is $1$). Moreover $f$ maps $[G(Q),G(Q)]$ isomorphically onto $2\cdot A_n$, which is perfect since its a stem extension of $A_n$ perfect ($n\geq 5$). Denote by $Ab:G(Q)\to \Z$ the map corresponding to the abelianisation. The map $f\times Ab :G(Q)\to 2\cdot A_n\times \Z$ is injective since $K\cap [G(Q),G(Q)]=1$. The image of $[G(Q),G(Q)]$ is $2\cdot A_n\times \{0\}$ and $Ab$ is surjective. Hence $f\times Ab$ is surjective. We have proved that $f\times Ab$ is an isomorphism.  
\end{proof}
\section{non connected $k$-cycle quandles}
We show that the second quandle homology of a finite injective quandle can be computed in a representation group of the quandle by giving a formula for that (first subsection). In the second subsection, we consider non connected subquandles of $Conj(S_n)$ consisting of $k$-cycles (with $n\neq 6,7$), we determine representation groups for these quandles, compute their second quandle homology using the formula and also determine their enveloping group.
\subsection{Quandle homology in representation groups}
Let $Q$ be a finite injective quandle. We have denoted by $G(Q)$ the enveloping group. Let $G(Q)_\alpha$ be a representaiton group of $Q$ and as in section $1$ denote by $\pi_\alpha$ the quotient map $G(Q)\to G(Q)_\alpha$. 
\begin{proposition}\label{pi}
Let $Q$ be a finite injective quandle :
\item[1)] $\pi_\alpha$ restricts to an isomorphism $[G(Q),G(Q)]\to [G(Q)_\alpha,G(Q)_\alpha]$.
\item[2)] For $x\in Q$, $C_{G(Q)}(e_x)=\pi_\alpha^{-1}( C_{G(Q)_\alpha}(\pi_\alpha(e_x)))$, where $C$ is for centralizers.
\item[3)] $\mathrm{ker}(\pi_\alpha) \simeq \Z^{n_Q}$, where $n_Q$ is the number of orbits of $Q$.
\end{proposition}
\begin{proof}
$1)$ was discussed in the second paragraph of subsection \ref{repgr}. We prove $2)$. We have seen in the second paragraph of subsection \ref{repgr} that $\mathrm{ker}(\pi_\alpha)\cap [G(Q),G(Q)]=1$. Take $g\in C_{G(Q)_\alpha}(\pi_\alpha(e_x))$ and let $\tilde{g}$ be a lift of $g$ to $G(Q)$. We have $[\tilde g,e_x] \in \mathrm{ker}(\pi_\alpha)\cap [G(Q),G(Q)]=1$. This proves that $\tilde g$ lies in $C_{G(Q)}(e_x)$. $2)$ follows since $G(Q)\to G(Q)_\alpha$ is a central extension. $3)$ was discussed in the first paragraph of subsection \ref{repgr}.
\end{proof}
Recall that the degree map $d:G(Q)\to \Z$ is the morphism given by $e_x\mapsto 1$ for $x\in Q$. We will use the notation $C_H(g)$ for the centralizer of $g$ in the group $H$. Set $A_x=C_{G(Q)}(e_x)\cap \mathrm{ker}(d)$.
\begin{proposition}
\item[1)] $\mathrm{ker}(\pi_\alpha)\cap \mathrm{ker}(d)\simeq \Z^{n_Q-1}$, where $n_Q$ is the number of orbits of $Q$.
\item[2)] $\mathrm{ker}(\pi_\alpha)\cap \mathrm{ker}(d)\simeq \Z^{n_Q-1}$ has finite index in $A_x$.
\item[3)] $\pi_\alpha([A_x,A_x])=[C_{G(Q)_\alpha}(\pi_\alpha(e_x)),C_{G(Q)_\alpha}(\pi_\alpha(e_x))]$.
\end{proposition}
\begin{proof}
Take $x\in Q$. $e_x^{a}$ for $a=\vert Inn(Q)\vert$ lies in the kernel of the morphism $G(Q)\to Inn(Q)$ given by $e_x\mapsto L_x$. Hence $e_x^a$ lies in the center $Z_Q$ of $G(Q)$ (Cf. section 1). Now $\ker(\pi_\alpha)\oplus Tor(Z_Q)=Z_Q$, where $Tor(Z_Q)$ is the torsion subgroup of $Z_Q$. Hence, $e_x^ab$ lies in $\ker(\pi_\alpha)$ for some $b\in Tor(Z_Q)$. $d(b)=0$ since $b$ is a torsion element and hence $d(e_x^ab)=a\neq 0$. This proves that $d$ is not trivial over $\ker(\pi_\alpha)$. Since $\mathrm{ker}(\pi_\alpha) \simeq \Z^{n_Q}$ (previous proposition) and $d$ is non trivial with image in $\Z$. $\mathrm{ker}(\pi_\alpha)\cap \mathrm{ker}(d)\simeq \Z^{n_Q-1}$. This proves $1)$. $2)$ follows from the fact that $\ker(\pi_\alpha)$ has finite index in $G(Q)$. It is clear that the image of $A_x$ under $\pi_\alpha$ lies in $C_{G(Q)_\alpha}(\pi_\alpha(e_x))$. Hence, $\pi_\alpha([A_x,A_x])\subset [C_{G(Q)_\alpha}(\pi_\alpha(e_x)),C_{G(Q)_\alpha}(\pi_\alpha(e_x))]$. Take $a$ and $b$ in $C_{G(Q)_\alpha}(\pi_\alpha(e_x))$. By $2)$ of the previous proposition $a$ and $b$ can be lifted to $\tilde{a}$ and $\tilde{b}$ in $C_{G(Q)}(e_x)$. Since $e_x$ is central in $C_{G(Q)}(e_x)$ :
$$ \pi_\alpha([\tilde a e_x^{-d(\tilde a)},\tilde b e_x^{-d(\tilde b)}])=\pi_\alpha([\tilde a, \tilde b])=[a,b],$$
where the bracket is for the commutators. This proves the inclusion :$$[C_{G(Q)_\alpha}(\pi_\alpha(e_x)),C_{G(Q)_\alpha}(\pi_\alpha(e_x))] \subset \pi_\alpha([A_x,A_x]).$$ Indeed, $d(\tilde a e_x^{-d(\tilde a)})=d (\tilde b e_x^{-d(\tilde b)})=0$ and both elements lie in $C_{G(Q)}(e_x)$. We have proved $3)$.
\end{proof}
\begin{proposition}
Let $Q$ be a finite injective quandle. For $x \in Q$ :
$$A_x^{ab}\simeq \Z^{n_Q-1} \bigoplus \frac{[G(Q)_\alpha,G(Q)_\alpha]\cap C_{G(Q)_\alpha}(\pi_\alpha(e_x))}{[C_{G(Q)_\alpha}(\pi_\alpha(e_x)),C_{G(Q)_\alpha}(\pi_\alpha(e_x))]},$$
where $n_Q$ is the number of orbits of $Q$.
\end{proposition}
\begin{proof}
Set $B_x=[G(Q),G(Q)] \cap A_x$. $B_x$ consist of the torsion elements of $A_x$ (propositon \ref{tor}). We have $[A_x,A_x] \subset B_x$ and hence we have an exact sequence :
$$1\to B_x/[A_x,A_x]\to A_x^{ab}\to A_x/B_x\to 1.$$
Notice that $A_x/B_x$ can be idetinfied to a subgroup of $G(Q)^{ab}=G(Q)/[G(Q),G(Q)]$, which is free abelian and hence $A_x/B_x$ is free abelian. This prove that the exact sequence above splits and that we have :
$$  A_x^{ab} \simeq   A_x/B_x  \oplus B_x/[A_x,A_x].$$
We have seen in the previous proposition that $\ker{\pi_\alpha} \cap \mathrm{ker}(d) \subset A_x$ is free abelian isomorphic to $\Z^{n_Q-1}$. $B_x$ is of torsion. Hence, $\ker(\pi_\alpha) \cap \mathrm{ker}(d)$ injects into $A_x/B_x$. Since (previous proposition) $\ker{\pi_\alpha}\cap \mathrm{ker}(d)$ has finite index in $A_x$ its image in $A_x/B_x$, the free abelian group, has finite index. This proves that $A_x/B_x \simeq \Z^{n_Q-1}$. Therfore : 
$$ A_x^{ab} \simeq \Z^{n_Q-1} \oplus B_x/[A_x,A_x].$$    
We will show that :
$$B_x/[A_x,A_x]\simeq \frac{[G(Q)_\alpha,G(Q)_\alpha]\cap C_{G(Q)_\alpha}(\pi_\alpha(e_x))}{[C_{G(Q)_\alpha}(\pi_\alpha(e_x)),C_{G(Q)_\alpha}(\pi_\alpha(e_x))]},$$
 proving the proposition. Using the inclusion $[G(Q),G(Q)] \subset \mathrm{ker}(d)$ then $2)$ of  propositon \ref{pi}, we get : 
$$B_x=[G(Q),G(Q)]\cap A_x=[G(Q),G(Q)]\cap C_{G(Q)}(e_x)=[G(Q),G(Q)] \cap \pi^{-1}(C_{G(Q)_\alpha}(\pi_\alpha(e_x))).$$
 By proposition \ref{pi}, $\pi_\alpha$ restricts to an isomorphism $[G(Q),G(Q)]\to [G(Q)_\alpha,G(Q)_\alpha]$. Hence, $\pi_\alpha$ maps $B_x$ isomorphically to $[G(Q)_\alpha,G(Q)_\alpha]\cap C_{G(Q)_\alpha}(\pi_\alpha(e_x))$. Now by $3)$ of the previous proposition 
$$\pi_\alpha([A_x,A_x])=[C_{G(Q)_\alpha}(\pi_\alpha(e_x)),C_{G(Q)_\alpha}(\pi_\alpha(e_x))].$$ 
Since $[A_x,A_x]\subset B_x$, the last two facts imply that :
$$B_x/[A_x,A_x]\simeq \frac{[G(Q)_\alpha,G(Q)_\alpha]\cap C_{G(Q)_\alpha}(\pi_\alpha(e_x))}{[C_{G(Q)_\alpha}(\pi_\alpha(e_x)),C_{G(Q)_\alpha}(\pi_\alpha(e_x))]},$$
We have proved the proposition.
\end{proof}

\begin{theorem}\label{Qh}
Let $(G,f)$ be a representation group for $Q$ a finite injective quandle. For each orbit in $Q$ pick an element $x_i$. The second quandle homology of $Q$ is given by :
$$H_2(Q)=\Z^{n_Q(n_Q-1)}\bigoplus (\bigoplus_i \frac{[G,G]\cap C_G(g_{x_i})}{[C_G(g_{x_i}),C_G(g_{x_i})]},$$
where $n_Q$ is the number of orbits of $Q$ and $g_{x_i}$ is any element in $f^{-1}(L_{x_i})$.
\end{theorem}
\begin{proof}
This follows from the last proposition and Eisermann's formula given in theorem \ref{eis0}. We note that different lifts of $L_{x_i}$ differ by central elements and hence define identical groups.
\end{proof}
\begin{remark}
The above formula for a finite connected injective quandle $Q$ gives $H_2(Q)=\frac{[G,G]\cap C_G(g_{x})}{[C_G(g_{x}),C_G(g_{x})]}$ and the formula in proposition \ref{RH} gives $([G,G]\cap C_G(g_{x}))^{ab}$. Hence they are isomorphic. We give a direct explanation. Let $\pi$ be the surjective projection from $G(Q)\to G$. $G$ is generated by the elements $\pi(e_x)$ for $x\in Q$ which are conjugate. Hence, $G^{ab}$ is cyclic generated by $\pi(e_x)$ for any $x$. We hence have $G=[G,G]\langle \pi(e_x) \rangle$ for any $x$. Since $\pi(e_x)$ is central in $C_G(\pi(e_x))$ we get from $G=[G,G]\langle \pi(e_x) \rangle$ that $$C_G(\pi(e_x))=([G,G]\cap C_G(\pi(e_x)))\langle \pi(e_x) \rangle$$ and hence $$[C_G(\pi(e_x)),C_G(\pi(e_x))]=[[G,G]\cap C_G(\pi(e_x)),[G,G]\cap C_G(\pi(e_x))].$$ This allows to identify both formulas.
\end{remark}
\subsection{Non connected $k$-cycle quandles}
Let $n\geq 3$ be an integer and take $k\in \{n,n-1\}$, $k$ odd. We will assume that $n\neq 6,7$.  We will denote by $Q$ the subquandle of $k$-cycles in $Conj(S_n)$. The quandles are not connected. They split exactly into two orbits under $Inn(Q)\simeq A_n$ for $n\geq 4$ and $Inn(Q)\simeq \{1\}$ for $n=3$. We will regard $Q$ as a subquandle of $A_n$ for $n\geq 4$. We recall that since it is a subquandle of a conjugacy quandle of a group, the quandle is injective. We will identify the map $Q\to Inn(Q),x\mapsto L_x$ to the inclusion $Q\to A_n$ for $n\geq 4$. 
\begin{theorem}
For $n\geq 4$ and $n\neq 6,7$, $2\cdot A_n$ is a representation group for $Q$ the subquanle of $k$-cycles of $Conj(S_n)$ for $k$ odd in $\{n,n-1\}$.
\end{theorem}
\begin{proof}
Arguing as in theorem \ref{x*}, we get that any element $x\in Q$ admits a unique lift $x^*$ of order $k$ and that we have a quandle morphism $Q\to \{x^* \vert x\in Q\},x\mapsto x^*$. We therefore have a morphism $G(Q)\to 2\cdot A_n$ given by $e_x\mapsto x^*$. It follows from proposition \ref{RS} that $2\cdot A_n$ is a representation group for $Q$ since in the range of $n$ in the theorem $2\cdot A_n$ is a Schur cover of $A_n\simeq Inn(Q)$.  
\end{proof}
\begin{theorem}
The second quandle homology of $Q$ the subquanle of $k$-cycles of $Conj(S_n)$ for $k$ odd in $\{n,n-1\}$ $(n\geq 3, n\neq 6,7)$ is given by:
$$ H_2(Q)\simeq \begin{cases} \Z^2 & \text{if } n=3\\
 \Z^2\bigoplus (\Z/2\Z)^2 & \text{if } n=4\\ 
\Z^2\bigoplus (\Z/2k\Z)^2 &\text{otherwise}
\end{cases}$$
\end{theorem}
\begin{proof}
For $n=3$, one can easly check that the enveloping group is $\Z^2$. Since we have $2$-orbit we get from Eisermann's formula in theorem \ref{eis0} that $H_2(Q)\simeq \ker(d)^2\simeq \Z^2$. We write as in the proof of the previous theorem $x^*$ for the unique lift of $x\in Q$ to $2\cdot A_n$ of order $k$. The centraliser $C_{2\cdot A_n}(x^*)$ lies in the preimage of the $C_{A_n}(x)$. Since $k\in\{n,n-1\}$, $C_{A_n}(x)$ is the cyclic group generated by $x$. Its preimage in $2\cdot A_n$ is the subgroup generated by $z$ and $x^*$. Hence, $C_{2\cdot A_n}(x^*)=\langle z,x^* \rangle$. This group is isomorphic to $\Z/2k\Z$. For $n\geq 5$, $A_n$ is perfect and hence $2\cdot A_n$ is perfect, we have $2$-orbits, and $\Z/2k\Z$ is abelian. Hence, by theorem \ref{Qh}, $H_2(Q)\simeq \Z^2\bigoplus (\Z/2k\Z)^2$ for $n\geq 5$. For $n=4$, $2\cdot A_n$ is not perfect. $[2\cdot A_n,2\cdot A_n]$ is the preimage of $[A_n,A_n]$ consisting of $1$ and three double transpositions. Hence $C_{2\cdot A_n}(x^*)\cap [2\cdot A_n,2\cdot A_n]=\langle z \rangle \simeq \Z/2\Z$. We have $2$-orbits and $C_{2\cdot A_n}(x^*)$ is abelian. We get using theorem \ref{Qh} that $H_2(Q)\simeq \Z^2\bigoplus (\Z/2\Z)^2$ for $n=4$.
\end{proof}
\begin{theorem}
For $n\geq 5$ and $n\neq 6,7$, the enveloping group of the $k$-cycle quandle of $S_n$ for $k$ odd in $\{n-1,n\}$ is isomorphic to $2\cdot A_n\times \Z^2$.
\end{theorem}
\begin{proof}
One argues as in the proof of theorem \ref{envCycle} but considering that we have $2$ orbits.
\end{proof}
\section{conjugacy class in $SL(n,q)$}
We prove that the subquandle consisting of a conjugacy class generating a Schur cover $G$ of a centerless group has representation group $G$. We use this to compute the second quandle homology of the subquandle of $Conj(SL(n,q))$ (with some exceptions for some $(n,q)$) corresponding to the conjugacy class of a diagonalizable element. The section is not divided into subsections.
\begin{lemma}
Let $G$ be a Schur cover of $H$. If $H$ is centerless, then $G$ is a Schur cover of $Inn(G)$.
\end{lemma}
\begin{proof}
This is immediate since $H\simeq Inn(G)$.
\end{proof}
\begin{proposition}
Let $G$ be a Schur cover of a centerless group $H$. Assume that the conjugacy class $Cl_G(x)$ of $x$ generates $G$. 
\item[1)] The subquandle $Q=Cl_G(x)$ of $Conj(G)$ is connected and $G$ is a representation group of $Cl_G(x)$.
\item[2)] The second quandle homology of $Cl_G(x)$ is isomorphic to $(C_G(x)\cap [G,G])^{ab}$. 
\end{proposition}
\begin{proof}
By the proposition \ref{In}, $Inn(Q) \simeq Inn(G) \simeq H$ and $Q$ is connected. The map $Q\to Inn(Q),x\mapsto L_x$ can be identified with the projection $Q\to H$. The quandle $Q$ is also injective by proposition \ref{inj}. We have a morphism $G(Q) \to G$ given by $e_x\mapsto x$ for $x\in Q$ and $x$ is a lift of the projection of $x$ into $H$. It follows from proposition \ref{RS} that $G$ is a representation group. This proves $1)$. $2)$ is an application of proposition \ref{RH}.
\end{proof}
 Take $n\geq 2$ and $q=p^m$ with $p$ a prime number and $m$ a positive integer. Recall that :
\begin{theorem}
$SL(n,q)$ is perfect for $(n,q)\neq (2,2),(2,3)$.
\end{theorem}
Recall that :
\begin{theorem}
For $(n,q) \neq (2,2),(2,3)$, $PSL(n,q)$ is simple.
\end{theorem}
The following result is due to the work of Steinberg and is covered in \cite{Sc}.
\begin{theorem}[Steinberg]
For $(n,q)\notin \{(2,4),(2,9),(3,2),(3,3),(3,4), (4,2) \}$, the group $SL(n,q)$ is a Schur cover of $PSL(n,q)$
\end{theorem}
We will say that $(n,q)$ is regular if $(n,q)$ does not figure in the exceptions of the last three theorems and $q\neq 2,3$. We will denote by $\F_q$ the field on $q$ elements and by $\F_q^*$ its group of invertibles.
\begin{proposition}
For $(n,q)$ regular, let $x$ be a diagonalizable element of $SL(n,q)$ and let $Q$ denotes the quandle corresponding to the conjugacy class of $x$. If $x$ is non scalar $SL(n,q)$ is a representation group for $Q$  and for any diagonalizable $x$ :
$$H_2(Q)\simeq (\Z/(q-1)\Z)^{m-1},$$
where $m$ is the number of eigenvalues of $x$. 
\end{proposition}
\begin{proof}
If $x$ is scalar, $Q=\{x\}$ and the quandle is connected. $G(Q)\simeq \Z$ and hence its derived group is trivial. It follows from Eisermann's formula for $H_2(Q)$ that $H_2(Q)=0$. Take $x$ not scalar, i.e. $m \geq 2$. Since $(n,q)$ is regular $SL(n,q)\to PSL(n,q)$ is a stem extension let $H\subset SL(n,q)$ be the subgroup generated by $Q$. $H$ is normal and non trivial and hence it maps onto $PSL(n,q)$ since $PSL(n,q)$ is simple. It follows from lemma \ref{Sub} that $H=SL(n,q)$. $PSL(n,q)$ is in particular centerless and the previous proposition applies to $Q$. We have :
$$H_2(Q)\simeq (C_{SL(n,q)}(x)\cap [SL(n,q),SL(n,q)])^{ab}.$$
Since $SL(n,q)$ is perfect for the regular values, $H_2(Q)\simeq C_{SL(n,q)}(x)^{ab}$. For $m$ eigenvalues the centraliser of $x$ in $GL(n,q)$ is given by a product $K=\prod_{i=1}^{ m} GL(n_i,q)$, where $\sum n_i=n$. And $C_{SL(n,q)(x)}$ is the kernel of $f:K\to \F_q^*$ given by $f(g_1,\dots,g_m)=det(g_1)\cdots det(g_m)$. Notice that $L=\prod_{i=1}^m SL(n_i,q)$ is a subgroup of $C_{SL(n,q)}(x)$. Since $q\neq 2,3$ (regular values), $SL(n_i,q)$ is perfect, $L= [L,L]\subset [C_{SL(n,q)}(x),C_{SL(n,q)}(x)]$ but $K/L\simeq (\F_q^*)^{m}$ ($q\neq 2$) is abelian and hence $C_{SL(n,q)}(x)^{ab}$ is the image of $C_{SL(n,q)}(x)$ in $K/L$.  $f:K\to \F_q^*$ factors through $K/L$ by the morphism $(\F_q^*)^m\to \F_q^*,(\lambda_1,\dots,\lambda_m)\mapsto \lambda_1\cdots\lambda_m$. Hence, $C_{SL(n,q)}(x)^{ab}$ is the subgroup of $(\lambda_1,\dots,\lambda_m)\in (\F_q^*)^m$ such that $\lambda_m^{-1}=\lambda_1\cdots \lambda_{m-1}$. This subgroup is isomorphic to $(\F_q^*)^{m-1}$ and $\F_q^*\simeq \Z/(q-1)\Z$. We have proved that :
$$ H_2(Q)\simeq C_{SL(n,q)}(x)^{ab} \simeq (\Z/(q-1)\Z)^{m-1}.$$
\end{proof}
\section{$Conj(G)$ for groups with trivial $B_0(G)$}
For $G$ a finite group with trivial Bogomolov multiplier, we give, in the first subsection, a formula to compute the second quandle homology of $Conj(G)$ in $G$. We apply this, in the second subsection, to compute the second quandle homology of $Conj(GL(2,q))$ with $q\neq 2,3$.  

\subsection{Groups with trivial Bogomolov multiplier}
In this section $G$ is a finite group. We will denote by $A(G)$ the enveloping group of $Conj(G)$. We recall that quandles $Conj(G)$ are injective.
Let $\pi :A(G)\to G$ be the group morphism given by $e_g\to g$. We will denote by $Z$ the kernel of $\pi$.
\begin{proposition}
The subgroup $Z$ is a finitely generated abelian group lying in the center of $A(G)$.
\end{proposition}
\begin{proof}
If we compose $\pi$ by the morphism $G\to Inn(G)$ mapping $g$ to the conjugacy by $g$ we get the central extension $A(G)\to Inn(Conj(G))=Inn(G)$ seen in the first section. This proves that $Z$ lies in the center of $A(G)$ and is abelian. Since $Z$ is of finite index in $A(G)$ and $A(G)$ is finitely presented, $Z$ is finitely presented and hence is finitely generated.
\end{proof}

The Bogomolov multiplier of a finite group $G$, denoted $B_0(G)$, is the subgroup of $H^2(G,\QZ)$ consisting of classes whose restrictions to $H^2(A,\QZ)$ are trivial for any abelian subgroup $A$ of $G$. Many classical groups has trivial trivial Bogolomov multiplier and (\cite{Bog}) simple groups and with some exceptions quasi-simple groups have trivial Bogomolov multiplier, i.e. $B_0(G)=0$.
\begin{proposition}\label{der}
If $B_0(G)=0$ then :
\item[1)] $Z\cap [A(G),A(G)]=1$.
\item[2)] $\pi$ restricts to an isomorphism $[A(G),A(G)]\to [G,G]$.
\item[3)] $Z\simeq \Z^{c(G)}$, where $c(G)$ is the number of conjugacy classes of $G$.
\end{proposition}
\begin{proof}
The intersection in $1)$ is the torsion subgroup of $Z$ (Cf. proposition \ref{tor}). This group is trivial if $B_0(G)=0$ (\cite{Irc}). $2)$ is an immediate consequence of $1)$. We have seen that the torsion of $Z$ is trivial if $B_0(G)=0$ so $Z$ is free abelian. Since the derived group of $A(G)$ is finite (Cf. proposition \ref{tor}), $Z\to A(G)^{ab}$ is injective. Now since $Z$ has finite index in $A(G)$ and $Z\cap [A(G),A(G)]=1$, the rank of $Z$ is the rank of $A(G)^{ab}$ and hence is equal to the number of orbits of $Conj(G)$. The number of orbits is the number of conjugacy classes of $G$. We have proved $3)$. 
\end{proof}
Recall that the degree map for $A(G)$ is the morphism $d:A(G)\to \Z$ given by $d(e_g)=1$, for $g\in G$.
\begin{proposition}\label{Zd}
If $B_0(G)=0$, then $Z\cap \mathrm{ker}(d)\simeq \Z^{c(G)-1}$ with $c(G)$ the number of conjugacy classes of $G$.
\end{proposition}
\begin{proof}
$e_g^{\vert G \vert}$ lies in the kernel of $\pi$ and $d(e_g^{\vert G \vert})={\vert G \vert}$. This proves that $d_{\vert Z}$ is not trivial. The proposition follows from this and the fact that $Z\simeq \Z^{c(G)}$ for $B_0(G)=0$.
\end{proof}
Set  $$A_g=C_{A(G)}(e_g)\cap \mathrm{ker}(d),$$
where $C_{A(G)}(e_g)$ is the centralizer of $e_g$ in $A(G)$.
\begin{proposition}
\item[1)] $Z\cap \mathrm{ker}(d)$ has finite index in $A_g$.
\item[2)] $\pi^{-1}(C_G(g))=C_{A(G)}(e_g)$ where $C$ is for centralizers in the corresponding groups.
\item[3)] The morphism $\pi$ restricts to a surjective morphism $A_g\to C_G(g)$, where $C_G(g)$ is the centralizer of $g \in G$.
\end{proposition}
\begin{proof}
We prove $1)$. $Z$ is a subgroup of $C_{A(G)}(e_g)$ and has finite index in $A(G)$ an hence finite index in the centralizer. $1)$ follows. The fact that $C_{A(G)}(e_g)$ map to $C_G(g)$ under $\pi$ is clear. We prove it is surjective. Take $h$ in $C_G(g)$. It follows from the definition of $A(G)$ that $e_h$ lies in $C_{A(G)}(e_g)$. This proves $2)$, since $A(G)\to G$ is a central extension. To prove $3)$ take $e_he_1^{-1}$. It lies in $C_{A(G)}(e_g)$ and in $\mathrm{ker}(d)$. Moreover $e_he_1^{-1}$ is mapped onto $h$ under $\pi$. We have proved $3)$. 
\end{proof}
\begin{proposition}
If $B_0(G)=0$, then the group $A_g^{ab}$ is isomorphic to :
$$\Z^{c(G)-1}\oplus \frac{[G,G]\cap C_G(g)}{[C_G(g),C_G(g)]},$$
where $c(G)$ is the number of conjugacy classes of $G$, $C_G(g)$ is the centralizer of $g$ in $G$.
\end{proposition}
\begin{proof}
Set $B_g=[A(G),A(G)] \cap A_g$. $B_g$ consist of the torsion elements of $A_g$ (propositon \ref{tor}). $[A_g,A_g] \subset B_g$ and hence we have an exact sequence.
$$1\to B_g/[A_g,A_g]\to A_g^{ab}\to A_g/B_g\to 1.$$
Notice that $A_g/B_g$ can be identified to a subgroup of $A(G)^{ab}=A(G)/[A(G),A(G)]$, which is free abelian and hence $A_g/B_g$ is free abelian. This prove that the exact sequence above splits and that we have :
$$  A_g^{ab} \simeq   A_g/B_g  \oplus B_g/[A_g,A_g].$$
We have (proposition \ref{Zd}) seen that $Z \cap \mathrm{ker}(d) \subset A_g$ is free abelian isomorphic to $\Z^{c(G)-1}$. $B_g$ is of torsion. Hence, $Z \cap \mathrm{ker}(d)$ injects into $A_g/B_g$. Since (previous proposition) $Z\cap \mathrm{ker}(d)$ has finite index in $A_g$ its image in $A_g/B_g$, the free abelian group, has finite index. This proves that $A_g/B_g \simeq \Z^{c(G)-1}$. Therefore : 
$$ A_g^{ab} \simeq \Z^{c(G)-1} \oplus B_g/[A_g,A_g].$$    
We will show that, $B_g/[A_g,A_g]\simeq \frac{[G,G]\cap C_G(g)}{[C_G(g),C_G(g)]}$, proving the proposition. Using the inclusion $[A(G),A(G)] \subset \mathrm{ker}(d)$ then $2)$ of the previous propositon, we get : 
$$B_g=[A(G),A(G)]\cap A_g=[A(G),A(G)]\cap C_{A(G)}(e_g)=[A(G),A(G)] \cap \pi^{-1}(C_G(g)).$$
 By proposition \ref{der} $\pi$ restricts to an isomorphism $[A(G),A(G)]\to [G,G]$. Hence, $\pi$ maps $B_g$ isomorphically to $[G,G]\cap C_G(g)$. Now by $3)$ of the previous proposition $\pi([A_g,A_g])=[C_G(g),C_G(g)]$. Since $[A_g,A_g]\subset B_g$, the last two facts imply that :
$$B_g/[A_g,A_g]\simeq \frac{[G,G]\cap C_G(g)}{[C_G(g),C_G(g)]}.$$
We have proved the proposition.
\end{proof}
\begin{theorem}\label{Gh}
Let $G$ be a finite group with trivial $B_0(G)$. The second quandle homology of $Conj(G)$ is given by :
$$H_2(Conj(G))\simeq \Z^{c(G)(c(G)-1)} \bigoplus (\bigoplus_{C\in Cl(G)}  \frac{[G,G]\cap C_G(g_C)}{[C_G(g_C),C_G(g_C)]}),$$
where $Cl(G)$ is the set of conjugacy classes of $G$, $g_C$ is a representative of the class $C$ and $C_G(g_C)$ is the centralizer of $g_C$ in $G$.
\end{theorem}
\begin{proof}
It follows immediately from the previous propositon and Eisermann's formula given in theorem \ref{eis0}.
\end{proof}
\begin{proposition}
Let $G$ be a finite perfect group with trivial Bogomolov multiplier. Using the notations of the previous theorem, the second quandle homology of $Conj(G)$ is given by :
$$H_2(Conj(G))\simeq \Z^{c(G)(c(G)-1)} \bigoplus (\bigoplus_{C\in Cl(G)}   C_G(g_C)^{ab})\simeq \Z^{c(G)(c(G)-1)} \bigoplus HH_1(\Z[G]) ,$$
where $HH_1(\Z[G])$ is the first Hochschild homology group of the group ring $\Z[G]$.
\end{proposition}
\begin{proof}
The first isomorphism in the equation is due to the previous theorem since $G$ is perfect. Now :
$$\bigoplus_{C\in Cl(G)}   C_G(g_C)^{ab}=\bigoplus_{C\in Cl(G)}   H_1(BC_G(g_C)),$$
where $BC_G(g_C)$ is the classifying space of the group $C_G(g_C)$. The last sum in the last equation is isomorphic to  $HH_1(\Z[G])$ (\cite{bur}). 
\end{proof}
\subsection{Conjugacy quandle of $GL(2,q)$}
In this section $q=p^n$ with $p$ a prime $n>0$ an integer. We will suppose that $q\neq 2,3$. We will denote by $\F_q$ the field with $q$ elements. 

\begin{lemma}
$PSL(2,q)$ and $PSL(3,4)$ are not isomorphic.
\end{lemma}
\begin{proof}
We have $\vert PSL(3,4) \vert=2^6\cdot 3^2\cdot 5\cdot 7$ and 
$$\vert PSL(2,q)\vert =\begin{cases} q(q^2-1) & \text{if } q \text{ is even}\\ \frac{q(q^2-1)}{2} & \text{ otherwise} \end{cases}.$$
Since $q$ and $(q^2-1)$ are coprime, the equality between cardinals of $PSL(3,4)$ and $PSL(2,q)$ for $q$ even implies that $q=2^6$. But $\vert PSL(2,2^6)\vert=2^6(2^{12}-1)=2^6(2^6-1)(2^6+1)=2^6\cdot 63 \cdot 65$ and $13$ divises $65$ but does not divide $\vert PSL(3,4) \vert$. This proves the lemma for $q$ even. For $q$ odd, since $q$ and $q^2-1$ are coprime the only possibilities for solutions of $$  \frac{q(q^2-1)}{2} =2^6\cdot 3^2\cdot 5\cdot 7$$ are $q=3^2$, $q=5$ and $q=7$. But 
$$ \vert PSL(2,3^2)\vert =2^3\cdot 3^2\cdot  5, \: \vert PSL(2,5) \vert =2^2\cdot 3 \cdot 5,\: \text{and } \vert PSL(2,7) \vert =2^3\cdot 3 \cdot 7.$$
All three values are different from $\vert PSL(3,4)\vert$. We have proved the lemma.

\end{proof}
The following statement should not be explained as $B_0(GL(2,q))\neq 0$ for $q=2,3$. We only give a proof for $q\neq 2,3$.
\begin{proposition}
For $q\neq 2, 3$, $B_0(GL(2,q))=0$.
\end{proposition}
\begin{proof}
Let $\alpha$ be a generator of $\F_q^*$. Consider the matrix :
$$A= \begin{pmatrix} \alpha & 0 \\ 0 & 1 \end{pmatrix}.$$
We have $GL(2,q)=SL(2,q)\rtimes H$, where $H=\langle A \rangle$ is the cyclic group generated by $A$. For $q\neq 2,3$, $SL(2,q)$ is perfect an hence (\cite{Bog}) we have :
$$ H^2(GL(2,q),\QZ)\simeq H^2(SL(2,q),\QZ)^H\times H^2(H,\QZ),$$
where the isomorphism is given by restriction maps. This implies that, for $q\neq 2,3$, $B_0(GL(2,q))$ injects in $B_0(SL(2,q))\times B_0(H)$. For $q\neq 2,3$, $SL(2,q)$ is quasi-simple (i.e.  $G$ perfect with $G/Z(G)$ simple) and $SL(2,q)$ is not a stem extension of $PSL(3,4)$. Indeed $PSL(3,4)$ is centerless and $PSL(2,q)$ is not isomorphic to $PSL(3,4)$ by the previous lemma. Hence, by the main theorem in \cite{Bog}, $B_0(SL(2,q))=0$. Now $H$ is cyclic and hence $B_0(H)=0$. Since both Bogomolov multipliers are trivial for $q\neq 2,3$, it follows from the inclusion of $B_0(GL(2,q))$ in $B_0(SL(2,q))\times B_0(H)$ that, for these values of $q$, $B_0(GL(2,q))=0$. This completes the proof of the proposition.  
\end{proof}
A conjugacy class of $GL(2,q)$ contains a matrix of exactly one of the following types :
\[ (1) \quad \begin{pmatrix} x & 0 \\ 0 & x \end{pmatrix} x\neq 0\quad (2) \quad \begin{pmatrix} x & 1 \\ 0 & x \end{pmatrix} \:x\neq 0 \quad (3)\quad \begin{pmatrix} x & 0\\ 0 & y \end{pmatrix} x\neq y, xy\neq 0. \]
\[ (4) \text{a matrix } M \text{ with irreducible characteristic polynomial} \]
Indeed the above type of matrices can be distinguished using the characteristic and minimal polynomials which are conjugacy invariants. We will say that the class is of type $(i)$ if it contains a matrix of type $(i)$. Recall that :
\begin{proposition}
We have $q-1$ classes of type $(1)$, $q-1$ classes of type $(2)$, $\frac{(q-1)(q-2)}{2}$ classes of type $(3)$ and $\frac{q(q-1)}{2}$ classes of type $(4)$. In total we have $q^2-1$ conjugacy classes in $GL(2,q)$.
\end{proposition}
Recall that :
\begin{proposition}
Let $A_i$ be a matrix of type $(i)$ and denote by $C(A_i)$ the centralizer of $A_i$ in $GL(2,q)$.
\item[1)] $C(A_1)=GL(2,q)$.
\item[2)] $C(A_2)=\{\begin{pmatrix} x & y \\ 0 & x \end{pmatrix} \vert x\neq 0 \}$.
\item[3)] $C(A_3)=\{\begin{pmatrix} x &0 \\ 0 & y \end{pmatrix}  xy\neq 0 \}$.
\item[4)] $C(A_4)=\F_q[A_4]^\times$.
\item[5)] The groups in $2),3),4)$ are abelian and hence have trivial derived subgroups.
\end{proposition}
\begin{proposition}
For $q\neq 2$, the derived subgroup $GL(2,q)'$ of $GL(2,q)$ is $SL(2,q)$.
\end{proposition}
\begin{proposition}
Let $A_i$ be a matrix of type $(i)$, denote by $C(A_i)$ the centralizer of $A_i$ in $GL(2,q)$ and by $GL(2,q)'$ derived subgroup of $GL(2,q)$. For $q\neq 2$ :
\item[1)] $C(A_1)\cap GL(2,q)'\simeq SL(2,q)$
\item[2)] $C(A_2)\cap GL(2,q)' \simeq \begin{cases}  \F_q & \text{if } q \text{ is even}\\\Z/2\Z \times \F_q & \text{otherwise}\end{cases}$
\item[3)] $C(A_3)\cap GL(2,q)'\simeq \Z/(q-1)\Z$.
\item[4)] $C(A_4)\cap GL(2,q)' \simeq \Z/(q+1) \Z$.
\end{proposition}
\begin{proof}
We only prove $4)$. We leave the others for the reader. Let $v$ be a non zero vector in $\F_q^2$. $(v,A_4v)$ forms a basis for $\F_q^2$ since $A_4$ has no eigenvalues. Let $f:\F_q[A_4]\to \F_q^2$ be the linear isomorphism given by $id\mapsto v$ and $A_4\mapsto A_4v$. Denote for $B\in F_q[A_4]$, $m_B$ the endomorphism of $F_q[A_4]$ corresponding to the multiplication by $B$. $fm_{A_4}f^{-1}$ correspond to $A_4$ and hence $fm_Bf^{-1}=B$ (we use the characteristic polynomial of $A_4$ to do the computation for $A_4$). It follows that $det(B)$ for $B\in \F_q[A_4]$ is equal to $det(m_B)$. Now $\F_q[A_4]\simeq \F_q[X]/(\chi_{A_4})$, where $\chi_{A_4}$ is the characteristique polynomial of $A_4$. $\chi_{A_4}$ is irreducible of degree $2$. Hence, $\F_q[A_4]/\F_q$ is a Galois extension of degree two with Galois group $\{1,\sigma\}$ where $\sigma(B)=B^q$. We have $det(m_B)=N_{\F_q[A_4]/\F_q}(B)=B\sigma(B)=B^{q+1}$ ( Here $\F_q$ is seen as $\F_q \cdot id$). Hence, for $B\in C(A_4)=\F_q[A_4]^\times$, $det(B)=B^{q+1}\in \F_q\cdot id$. It follows that $C(A_4)\cap GL(2,q)'$ is the set of $B$ such that $B^{q+1}=id$. Now $\F_q[A_4]^\times\simeq (F_q[X]/(\chi_{A_4}))^\times$ is cyclic of order $q^2-1$ and hence the equation $B^{q+1}=id$ defines a cyclic subgroup of order $q+1$. This proves $4)$. 
\end{proof}
Combining what we have seen in this subsection with theorem \ref{Gh} we get :

\begin{theorem}
\item[1)] For $q$ even and $q\neq 2$ :
$$H_2(Conj(GL(2,q)))\simeq \Z^{(q^2-1)(q^2-2)}\oplus \F_q^{q-1}\oplus (\Z/(q-1)\Z)^{\frac{(q-1)(q-2)}{2}}\oplus (\Z/(q+1) \Z)^{\frac{q(q-1)}{2}}.$$
\item[2)] For $q$ odd, $q\neq 3$ :
$$ H_2(Conj(GL(2,q)))\simeq \Z^{(q^2-1)(q^2-2)}\oplus (\Z/2\Z\times \F_q)^{q-1}\oplus (\Z/(q-1)\Z)^{\frac{(q-1)(q-2)}{2}}\oplus (\Z/(q+1) \Z)^{\frac{q(q-1)}{2}} $$

\end{theorem}
The last proposition allows one to determine the centralizers in $SL(2,q)$. The count and type of conjugacy class for $SL(2,q)$ can be done as for $GL(2,q)$. Hence, one can easily compute $H_2(Conj(SL(2,q)))$ for general values of $q$ as we have done for $GL(2,q)$.

 \end{document}